\documentclass[times]{article}

\usepackage{amsmath}
\usepackage{url}
\usepackage{cite}
\usepackage[amssymb]{SIunits}
\usepackage{lscape}
\usepackage{graphicx}

\newcommand{\D}{\mathrm{d}}  

\begin{document}

\bibliographystyle{unsrt}



\title{Quadrature for second-order triangles in the Boundary Element
  Method}
\author{Michael Carley}






\maketitle

\begin{abstract}
  A quadrature method for second-order, curved triangular elements in
  the Boundary Element Method (BEM) is presented, based on a polar
  coordinate transformation, combined with elementary geometric
  operations. The numerical performance of the method is presented
  using results from solution of the Laplace equation on a cat's eye
  geometry which show an error of order $P^{-1.6}$, where $P$ is the
  number of elements. 
\end{abstract}


\vspace{-6pt}

\section{INTRODUCTION}
\label{sec:intro}
\vspace{-2pt}

The Boundary Element Method (BEM), often called the `panel method' in
fluid dynamics, is a standard technique for the solution of boundary
integral equations in a number of fields. Historically, relatively
low-order discretizations have been used with geometries modelled
using first order elements, and surface variables modelled to zero or
first order on those elements. There are many methods available for
the computation of potential integrals on linear panels~\cite[for
example]{okon-harrington82b,okon-harrington82a,newman86,suh00,salvadori10,%
  carley13}, and their behaviour is reasonably well-understood,
allowing them to be implemented with confidence in production codes.

More recently, however, there has been increasing interest in the use
of higher order methods, in part to achieve better geometric fidelity,
and in part to improve the modelling of solutions. For example, a
recently developed panel method for whole aircraft
aerodynamics~\cite{willis-peraire-white05,willis-peraire-white07},
employing accelerated integration and summation techniques, depends on
the availability of a robust integration scheme for second order
panels, developed by the authors~\cite{willis-peraire-white06}, to
avoid some of the deficiencies inherent in other curved panel
integration techniques~\cite[for example]{wang-newman-white00}.

To clarify the application, we consider the solution of a boundary
integral formulation of the Laplace equation:
\begin{align}
  \label{equ:potential}
  \phi(\mathbf{x}) &= 
  \int_{S} 
  \frac{\partial\phi(\mathbf{y})}{\partial n}G(\mathbf{x},\mathbf{y})
  -
  \frac{\partial G(\mathbf{x},\mathbf{y})}{\partial n}\phi(\mathbf{y})
  \,\D S,
\end{align}
where $\phi$ denotes potential, $\mathbf{x}$ field point position,
$\mathbf{y}$ position on the surface $S$, and $n$ the outward pointing
normal to the surface. The Green's function $G$ is:
\begin{align}
  \label{equ:laplace}
  G(\mathbf{x};\,\mathbf{y}) &= \frac{1}{4\pi R},\\
  R &= |\mathbf{x}-\mathbf{y}|. \nonumber
\end{align}
To solve this problem using a BEM, the surface $S$ is discretized into
a number of elements, triangular in this case, over which $\phi$ is
approximated by some interpolant. This results in a linear system:
\begin{align}
  \label{equ:potential:system}
  c_{i}(\mathbf{x})\phi_{i} &= 
  \sum_{j=1}^{P}
  \int_{S_{j}}
  \frac{\partial\phi(\mathbf{y})}{\partial n}G(\mathbf{x},\mathbf{y})
  \,\D S_{j}
  -
  \sum_{j=1}^{P}
  \int_{S_{j}}
  \frac{\partial G(\mathbf{x},\mathbf{y})}{\partial n}\phi(\mathbf{y})  
  \,\D S_{j},
\end{align}
where $P$ is the number of elements (panels), $i$ is the index of a
surface point, $S_{j}$ is the surface of panel $j$, and the constant
$c_{i}(\mathbf{x})$ is a geometric property given by:
\begin{align}
  \label{equ:constant}
  c_{i} &=
  1 + 
  \int_{S} \frac{\partial G(\mathbf{x}_{i},\mathbf{y)}}{\partial n}\,\D S, 
\end{align}
equal to $1/2$ at a smooth point on the surface, taking some other
value at sharp edges. Inserting the interpolant for each element into
Equation~\ref{equ:potential:system} yields a system of equations
relating $\phi$ to $\partial\phi/\partial n$, allowing the problem to
be solved subject to the specification of some boundary condition. In
aerodynamic problems, this will usually be the Neumann boundary
condition, specifying the surface normal velocity
$\partial\phi/\partial n$. Upon solving for surface potential $\phi$,
the boundary integral can then be used to compute the potential or its
derivatives, i.e.\ fluid velocity, external to the surface, using
Equation~\ref{equ:potential}.

The core of the implementation is then the evaluation of the panel
integrals in Equation~\ref{equ:potential:system}. For planar panels,
there is no great difficulty, and for the Laplace equation, the
integration can be performed analytically using a variety of
approaches~\cite{okon-harrington82b,okon-harrington82a,newman86,%
  suh00,salvadori10,carley13}, although a numerical method is still
necessary for the Helmholtz equation for acoustic scattering and
radiation. When the panel is curved, however, a fully numerical method
is required, and extra difficulties arise in finding a transformation
which maps the curved panel to a reference domain where standard
quadrature rules can be applied. 

A common approach in integrating over planar panels is to convert to
polar coordinates with axis perpendicular to the element plane, which
mitigates difficulties caused by the $1/R$ singularity in the Green's
function. Such an approach has been used in computing the self-term
for curved panels~\cite{willis-peraire-white06}, and a similar
technique will be used here. Alternatives which have been used include
the mapping of the element onto a plane
triangle~\cite{wang-newman-white00}, which can, however, only be used
for the single layer potential, and onto a
sphere~\cite{willis-peraire-white06}, with appropriate conversions
between the appropriate Jacobians. 

In this paper, we present a technique for the evaluation of a
quadrature rule for second order curved panels which uses a polar
transformation of the integral combined with basic geometric
operations, to give a method more akin to the current techniques for
planar panels, but with additional complexity due to the need to
perform the geometric operations on curved element edges.

\section{ANALYSIS}
\label{sec:analysis}
\vspace{-2pt}

The quadrature method for the curved element is derived for a panel in
a reference position. For a planar element, there is no difficulty in
defining an element plane which can be used to fix a coordinate system
for the triangle, but for the curved elements we consider here, there
is clearly a choice to be made. A standard approach is to use a plane
tangent to the element, through some appropriate point, for example,
the field point $\mathbf{x}$ when a self-term is being computed, but
here we use the plane defined by the three corners of the triangular
element. The problem axes are rotated and shifted so that the corners
lie in a plane $z=0$ and the field point $\mathbf{x}=(0,0,z)$, i.e.\
the origin of the coordinate system is taken as the projection of the
field point onto the triangle reference plane. In this orientation,
cylindrical polar coordinates can be readily defined and used to carry
out the required integration. The rest of this section describes the
geometric operations employed, and the technique used to define the
quadrature rule.

\subsection{Description of second order triangles}
\label{sec:triangles}

\begin{figure}
  \centering
  \begin{tabular}{cc}
    \includegraphics{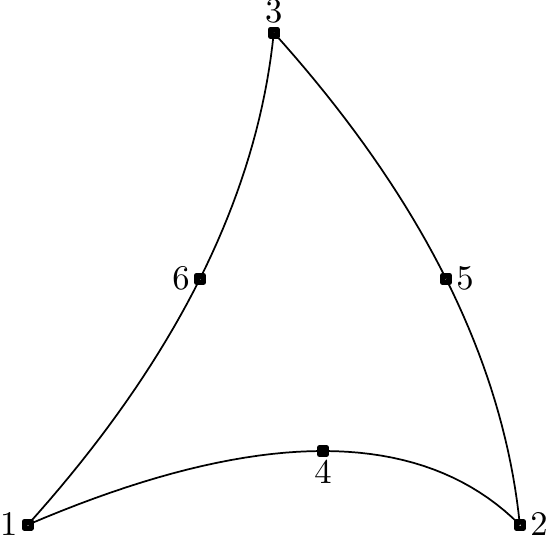} &
    \includegraphics{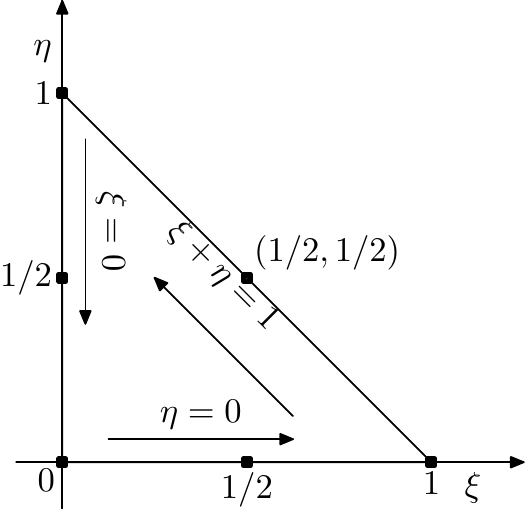}
  \end{tabular}
  \caption{Description of second-order curved triangle. Left hand
    image: curved triangle; right hand image: reference right-angle
    triangle.}
  \label{fig:triangle}
\end{figure}

Figure~\ref{fig:triangle} shows the notation used for description of
the second order triangle. Nodes are numbered~1,2,3 for the corners,
and~4,5,6 for the points internal to an edge. Quantities on the
element, including position, are interpolated using the second order
shape functions for the reference element:
\begin{align}
  \label{equ:interp}
  \mathbf{y} &= \sum_{i=1}^{6} L_{i}(\xi,\eta) \mathbf{y}_{i},\\
  \phi &= \sum_{i=1}^{6} L_{i}(\xi,\eta) \phi_{i},
\end{align}
where:
\begin{subequations}
  \label{equ:shape}
  \begin{align}
    L_{1} &= 2(1-\xi-\eta)(1/2-\xi-\eta),\\
    L_{2} &= 2\xi(\xi-1/2),\\
    L_{3} &= 2\eta(\eta-1/2),\\
    L_{4} &= 4\xi(1-\xi-\eta),\\
    L_{5} &= 4\xi\eta,\\
    L_{6} &= 4\eta(1-\xi-\eta),\\
    0&\leq\xi\leq 1,\,0\leq\eta\leq 1-\xi.\nonumber
  \end{align}
\end{subequations}
As shown in Figure~\ref{fig:triangle}, the edges of the triangle are
defined by second order interpolation on three points, and can be
described using a single variable $\gamma$, $0\leq\gamma\leq1$, with
$\gamma$ increasing in the anti-clockwise direction on each
edge. Inserting the conditions for each edge shown in
Figure~\ref{fig:triangle} gives three shape functions:
\begin{subequations}
  \begin{align}
    J_{1} &= 2\gamma^{2} - 3\gamma + 1,\\
    J_{2} &= 2\gamma^{2} - \gamma,\\
    J_{3} &= -4\gamma^{2} + 4\gamma,
  \end{align}
\end{subequations}
with the edge described by:
\begin{align}
  \label{equ:edge}
  \mathbf{y}(\gamma) &= 
  \mathbf{y}_{i}J_{1}(\gamma) +
  \mathbf{y}_{j}J_{2}(\gamma) +
  \mathbf{y}_{i+3}J_{3}(\gamma),
\end{align}
where $(i,j)$=$(1,2)$, $(2,3)$, $(3,1)$ for each edge respectively. A
point on an edge can be given in the general coordinates $(\xi,\eta)$
using the relations:
\begin{align}
  \label{equ:edge:area}
  (\xi,\eta) &=
  \left\{
    \begin{matrix}
      &(\gamma,0)\quad &&i=1;\\
      &(1-\gamma,\gamma),\quad &&i=2;\\
      &(0,1-\gamma),\quad &&i=3.
    \end{matrix}
  \right.
\end{align}
These shape functions will prove useful in determining intersections
between edges and lines in the plane, necessary in finding the domain
of integration for quadrature over the triangular element. 

\subsection{Integration over the triangle}
\label{sec:integration}

\begin{figure}
  \centering
  \includegraphics{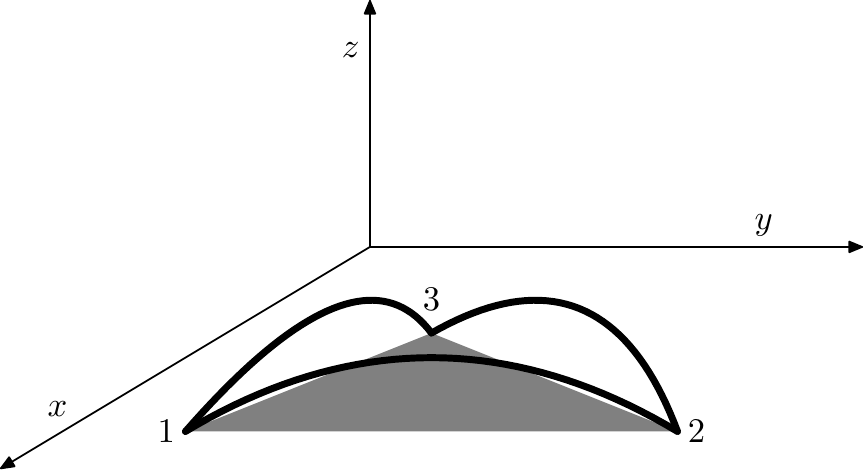}
  \caption{Coordinate system for calculation of integrals. The
    triangle is rotated so that the vertices 1,2,3 lie in the plane
    $z=0$ (shown shaded), and translated so that the field point lies
    at a position $(0,0,z)$.}
  \label{fig:orientation}
\end{figure}

Figure~\ref{fig:orientation} shows the curved triangle in its
reference position, rotated so that the corners lie in the plane
$z=0$, and shifted so that the field point lies at
$\mathbf{x}=(0,0,z)$. The integral to be evaluated is
\begin{align}
  \label{equ:potential:1}
  I &= \int_{0}^{1}\int_{0}^{1-\xi} f(\xi,\eta)
  J(\xi,\eta)\,\D\eta\,\D\xi,
\end{align}
where $J(\xi,\eta)$ is the Jacobian for the transformation from
$(\xi,\eta)$ to coordinates on the element surface. Converting to
Cartesian coordinates in the problem system of axes:
\begin{align}
  \label{equ:potential:2}
  I &= \int_{\triangle} f(\xi,\eta) r\,\D r\,\D\theta,
\end{align}
where integration takes place over the projection of the curved
triangle into the plane $z=0$ and the function $f(\xi,\eta)$ is
computed by transformation from $(r,\theta)$ to $(\xi,\eta)$.

\begin{figure}
  \centering
  \includegraphics{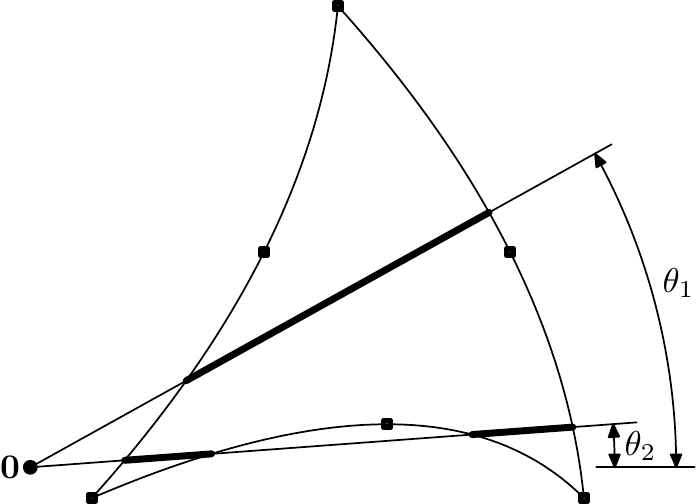}
  \caption{Intersection of rays with curved triangle: bold lines show
    the part of each ray which lies inside the domain of integration.}
  \label{fig:intersection}
\end{figure}

The integration of Equation~\ref{equ:potential:2} is conceptually
simple, and is readily applied to first order
elements~\cite{carley13}, but gives rise to some extra complexities in
the second order case, shown in Figure~\ref{fig:intersection}. As
written, the integration is composed of a sequence of integrals over
$r$, along rays at fixed values of $\theta$. In
Figure~\ref{fig:intersection}, two such rays are shown, at angles
$\theta_{1}$ and $\theta_{2}$. The ray at $\theta_{1}$ presents no
particular difficulties: it has one entry and one exit point on the
element boundary, both easily found using analytical methods (see
Section~\ref{sec:operations}). The ray at $\theta_{2}$, however, is
broken as it crosses the triangle boundary, having two entry and two
exit points. In evaluating Equation~\ref{equ:potential:2}, this case
must be handled, as must the case of a ray which lies tangent to an
edge. The algorithm of Section~\ref{sec:algorithm} handles these
special cases, using the geometrical operations of the next section. 

\subsection{Geometrical operations}
\label{sec:operations}

The quadrature algorithm, presented in Section~\ref{sec:algorithm},
depends on the availability of a number of elementary geometrical
operations, described in this section. These operations can be
implemented analytically using standard methods and are used to
determine the limits of integration in $\theta$, to break the integral
at possible points of discontinuity, and in $r$, to find the entry and
exit points on the triangle boundary.

The first operation is finding the intersection, $\gamma$ or
$(r,\theta)$, between a ray through the origin of angle $\theta$ and
an edge of the triangle. The coordinates of a point on the edge are
given by Equation~\ref{equ:edge}:
\begin{align}
  \label{equ:edge:1}
  x =
  x_{i}J_{1}(\gamma) + x_{j}J_{2}(\gamma) + x_{i+3}J_{3}(\gamma),\,
  y = y_{i}J_{1}(\gamma) + y_{j}J_{2}(\gamma) + y_{i+3}J_{3}(\gamma),
\end{align}
while the ray is given by $x=r\cos\theta$, $y=r\sin\theta$, so that,
upon substitution:
\begin{align}
  \label{equ:edge:2}
  (x_{i}\sin\theta - y_{i}\cos\theta)J_{1}(\gamma) 
  +
  (x_{j}\sin\theta - y_{j}\cos\theta)J_{2}(\gamma) 
  +
  (x_{i+3}\sin\theta - y_{i+3}\cos\theta)J_{3}(\gamma) 
  = 
  0,
\end{align}
which is a quadratic in $\gamma$. To find the intersection:
\begin{enumerate}
\item solve Equation~\ref{equ:edge:2} for $\gamma$;
\item for each value of $\gamma$, with $0\leq\gamma\leq1$
  \begin{enumerate}
  \item compute $(x,y)$ and $r=x/\cos\theta$ or $r=y/\sin\theta$;
  \item if $r>0$, $\gamma$ is a valid intersection.
  \end{enumerate}
\end{enumerate}
In this operation, $r$ is computed as shown in order to accept only
$r>0$, to avoid double counting of intersections. The condition
$0\leq\gamma\leq1$ is imposed to exclude corners of the triangle, as
these are handled separately in the algorithm. 

Tangents which pass through the origin must also be determined, in
order to break the integration at these points. They are found using
the equation of a tangent to a point $(x_{0},y_{0})$ on the edge:
\begin{align}
  \label{equ:edge:3}
  y &= y_{0} + 
  \left.\frac{\D y}{\D x}\right|_{x=x_{0}}(x-x_{0}),\\
  \left.
  \frac{\D y}{\D x}
  \right|_{x=x_{0}}
  &= 
  \frac{y_{i}J'_{1}(\gamma) + y_{j}J'_{2}(\gamma) +
    y_{i+3}J'_{3}(\gamma)}
  {x_{i}J'_{1}(\gamma) + x_{j}J'_{2}(\gamma) + x_{i+3}J'_{3}(\gamma)},
\end{align}
where the prime denotes differentiation with respect to
$\gamma$. Setting $x=y=0$ to find a tangent through the origin yields:
\begin{align}
  &x_{i}y_{j}
  \left[
    J_{1}'(\gamma)J_{2}(\gamma) -
    J_{2}'(\gamma)J_{1}(\gamma)
  \right] 
  +
  x_{j}y_{i+3}
  \left[
    J_{2}'(\gamma)J_{3}(\gamma) -
    J_{3}'(\gamma)J_{2}(\gamma)
  \right] \nonumber\\
  \label{equ:edge:5}
  &+
  x_{i+3}y_{i}
  \left[
    J_{3}'(\gamma)J_{1}(\gamma) -
    J_{1}'(\gamma)J_{3}(\gamma)
  \right]   
  = 0,
\end{align}
which is a cubic which can be solved for $\gamma$ subject to the
constraint $0<\gamma<1$ and that $\gamma$ be real. 

The final part of determining the limits of integration in $\theta$ is
the angle $\psi$ of a tangent to an edge, given by:
\begin{align}
  \label{equ:edge:6}
  \psi &= \tan^{-1}
  \frac{y_{i}J'_{1}(\gamma) + y_{j}J'_{2}(\gamma) +
    y_{i+3}J'_{3}(\gamma)}
  {x_{i}J'_{1}(\gamma) + x_{j}J'_{2}(\gamma) + x_{i+3}J'_{3}(\gamma)},
\end{align}
from the slope of the curve at a point $\gamma$ on an edge. A second
angle $\psi+\pi$ is also included in order to ensure that rays in both
tangent directions are included. 

\subsection{Quadrature algorithm}
\label{sec:algorithm}

The algorithm for the quadrature rule consists of a first stage in
which the range of integration in $\theta$ is broken into a set of
intervals, and a second in which the integration is performed over
these intervals. Initially it must be determined whether the origin
lies inside, outside or on the boundary of the triangle projected into
the plane $z=0$. This is done by first checking if the origin lies
outside a box containing the points $(x_{1},y_{1})$, $(x_{2},y_{2})$,
and~$(x_{3},y_{3})$. If not, its coordinates $(\xi,\eta)$ in the
triangle are found using Newton's method, and it is checked whether
they lie within the triangle.

In the first stage of processing:
\begin{enumerate}
\item find possible limits of integration $\theta_{i}$,
  $i=1,\ldots,N$, as the angles of rays joining the origin to the
  triangle's corners, the angles of tangents through the origin, and,
  if the origin lies on an edge, the angles of tangents to the edge;
\item adjust all angles to lie in the range $0,2\pi$;
\item sort the list of limits $\theta_{i}$ in ascending order;
\item if the origin lies inside the triangle, append the angle
  $\theta_{1}+2\pi$.
\end{enumerate}

Given the list of $N$ angles from the first stage, the nodes and
weights of the quadrature rule are found for each pair of limits,
$\theta_{i}$, $\theta_{i+1}$ by this procedure:
\begin{enumerate}
\item select a quadrature rule with abscissae $\theta_{k}$ and weights
  $w_{k}^{\theta}$, $k=1,\ldots,K$;
\item for each $k$:
  \begin{enumerate}
  \item find the radii $r_{j}$ of the intersections of a ray of angle
    $\theta_{k}$ with the triangle edges (prepend $r=0$ if the origin
    lies inside the triangle);
  \item for each pair of limits $r_{2j}$ and $r_{2j+1}$, select a
    quadrature rule $(r_{m},w_{m}^{r})$, $m=1,\ldots,M$;
  \item for each point $(r_{m}\cos\theta_{k},r_{m}\sin\theta_{k})$:
    \begin{enumerate}
    \item find the corresponding coordinates $(\xi,\eta)$ using
      Newton's method;
    \item append to the quadrature rule the abscissa $(\xi,\eta)$ and
      the weight $r_{m}w_{m}^{r}w_{k}^{\theta}/J_{2}(\xi,\eta)$, where
      $J_{2}$ is the Jacobian for conversion from $(\xi,\eta)$ to
      $(r,\theta)$. 
    \end{enumerate}
  \end{enumerate}
\end{enumerate}
The resulting quadrature rule can be used to evaluate an integral on
the panel by summation:
\begin{align}
  \label{equ:rule}
  \int_{0}^{1}\int_{0}^{1-\xi} f(\xi,\eta) J(\xi,\eta)\,\D\eta\,\D\xi
  &\approx \sum_{n} f(\xi_{n},\eta_{n}) J(\xi_{n},\eta_{n})w_{n},
\end{align}
where $J(\xi,\eta)$ is the Jacobian for conversion from $(\xi,\eta)$
to the surface coordinates $(x,y,z)$, allowing the rule to be used in
the same manner as standard quadratures for triangles. We note that
since the quadrature rule is mapped to the reference triangle, the sum
of the weights $w_{n}$ should be equal to $1/2$, which gives a
convenient error measure for checking the accuracy of the quadrature. 

Finally, a good starting guess for Newton's method in finding
coordinates on the reference triangle is provided by the intersection
point $(r,\theta)$ since its value of $\gamma$ on the edge is known,
and can be converted to $(\xi,\eta)$ using
Equation~\ref{equ:edge:area}. Each set of coordinates $(\xi,\eta)$ on
a ray can then be used as an initial guess for the evaluation at the
next quadrature point. 

\subsection{Quadrature selection}
\label{sec:quadrature}

The quadrature rule is implemented using a sequence of Gaussian
quadratures for integration in $\theta$ and $r$. This has the first
advantage that the endpoints of the integral are not included. Since
tangents are used to fix limits of integration, there is then no
ambiguity in determining the number of entry and exit points in
radius.

The quadrature rules are selected using a criterion which gives some
adaptivity to the intervals of integration. A point separation
parameter $\Delta\theta$ is specified and used to estimate the number
of points required in the quadrature rule:
\begin{align}
  \label{equ:selection}
  K &= \left[
    (\theta_{i+1}-\theta_{i})/\Delta\theta
  \right],
\end{align}
where $[\cdot]$ denotes rounding of the value to the nearest natural
number. The resulting value $K$ is then adjusted to lie in a range
$K_{\min}\leq K \leq K_{\max}$. A similar method is used to select
quadrature rules in $r$. In the calculations presented here, the
values of $\Delta\theta$ and $\Delta r$ are computed with user-defined
constants $N_{\theta}$ and $N_{r}$, and set as follows:
\begin{subequations}
  \begin{align}
    \Delta\theta &= \frac{\max_{i}\phi_{i}}{N_{\theta}},\\
    \Delta r &= \frac{\max_{i}\ell_{i}}{N_{r}},
  \end{align}
\end{subequations}
where $\phi_{i}$ is the angle subtended by the corner of the triangle
at vertex $i$, $i=1,2,3$ and $\ell_{i}$ is the length of the straight
edge starting at corner $i$. This gives a quickly computed abscissa
density under user control. 

Finally, in using the algorithm in a code, a criterion is required in
deciding when to use it, and when not. In this case, a parameter
$\sigma$ is computed based on easily evaluated geometric properties of
the element. These are the mean values of the nodes and the radius of a
sphere containing the element:
\begin{subequations}
  \begin{align}
    \overline{\mathbf{y}} &= \frac{1}{6}\sum_{i=1}^{6}\mathbf{y}_{i},\\
    \rho &= s\max_{i}|\mathbf{y}_{i}-\overline{\mathbf{y}}|,
  \end{align}
\end{subequations}
where $s$ is a scaling factor, with, in this case, $s=2^{1/2}$.  Given
these values for the element, $\sigma$ is determined as follows:
\begin{enumerate}
\item if $\mathbf{x}$ lies on the element, $\sigma=0$;
\item compute $\rho_{x}=|\mathbf{x}-\overline{\mathbf{y}}|$: if
  $\rho_{x}>\rho$, $\sigma=\rho_{x}/\rho$;
\item otherwise, $\sigma=z/\rho$,
\end{enumerate}
where $z$ is the distance of $\mathbf{x}$ from the element reference
plane, as noted above. This gives a quickly-computed parameter which
varies from~$0$ on the element to large values away from the element,
but remains small in some reasonable neighbourhood, so that it can be
used to select quadrature methods.

\section{NUMERICAL TESTS}
\label{sec:tests}
\vspace{-2pt}


The algorithm of the previous section is demonstrated using two sets
of results. The first is an illustration of the distribution of
quadrature nodes, using a sample element, while the second is an
assessment of the accuracy and convergence of the method when
implemented in a BEM program.

\subsection{Quadrature points}
\label{sec:points}

\begin{figure}
  \centering
  \begin{tabular}{ccc}
    \includegraphics{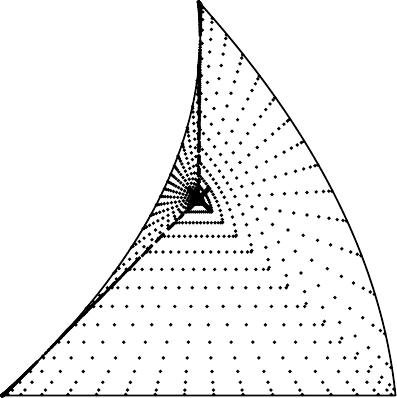} &
    \raisebox{-1mm}{\includegraphics{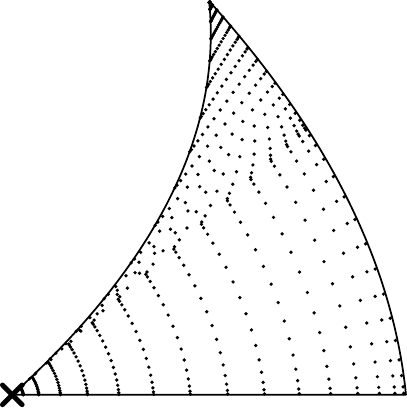}} &
    \includegraphics{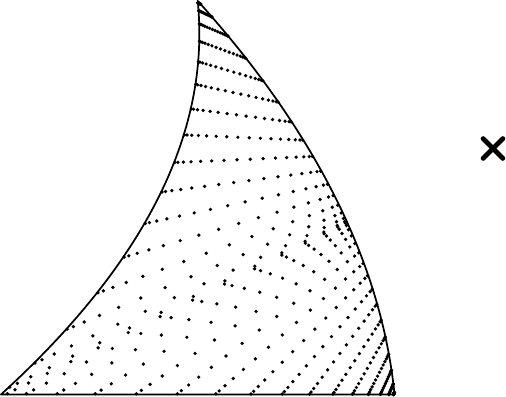} \\
    \raisebox{-1mm}{\includegraphics{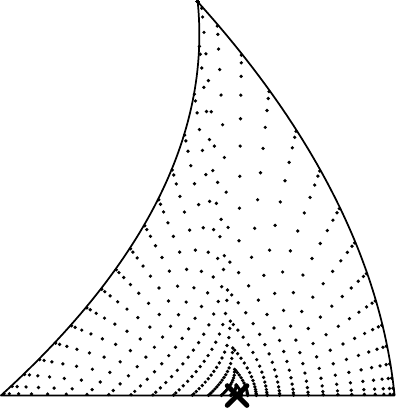}} &
    \includegraphics{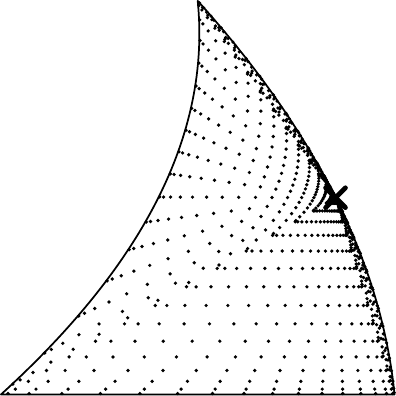} &
    \includegraphics{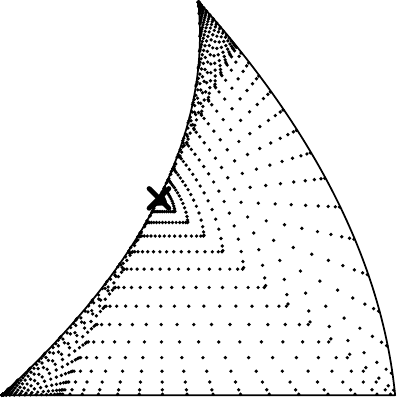} 
  \end{tabular}
  \caption{Quadrature points for different field points: origin is
    shown as a cross, quadrature points as dots. Top row: origin
    inside element, on vertex, and outside element. Bottom row: origin
    on straight, convex and concave edge.}
  \label{fig:points}
\end{figure}

To demonstrate the nature of the quadrature point distributions
generated by the algorithm, an element with one straight, one convex
and one concave edge has been used, with the origin placed inside and
outside the element, on a corner, and on each of the edges in
turn. In order to show the point distribution more clearly, fixed
length quadrature rules have been used, with sixteen points in both
angle and radius. This gives a higher than normal density in small
intervals of integration, and a lower density in larger regions, but
is helpful for visualization.

Figure~\ref{fig:points} shows the resulting quadratures, with the
origin indicated by a cross. Each of the cases considered gives rise
to a qualitatively different point distribution. With the origin
located inside the element, the region of integration is divided by
rays to each of the corners (compare Figure~8 of
reference~\cite{willis-peraire-white06}). When the origin is placed on
a corner of the element, there are two clearly demarcated domains of
integration, separated by the tangent to the curved edge at that
corner. Similarly, when the origin is moved outside the element, there
are two domains, separated in this case by a ray joining the origin to
the most distant corner. 

When the origin lies on an edge, the situation is slightly more
complicated. When it is on the straight edge, the element is divided
into two regions, separated by the ray to the furthest corner. When it
lies on the convex edge, there is a thin region of integration bounded
by the edge and the rays to the vertices on that edge. Conversely, on
the concave edge, the narrow region of integration is bounded by the
edge and the tangent to the edge at the origin.

\subsection{Numerical accuracy and convergence}
\label{sec:accuracy}

\begin{figure}
  \centering
  \includegraphics[width=0.6\textwidth]{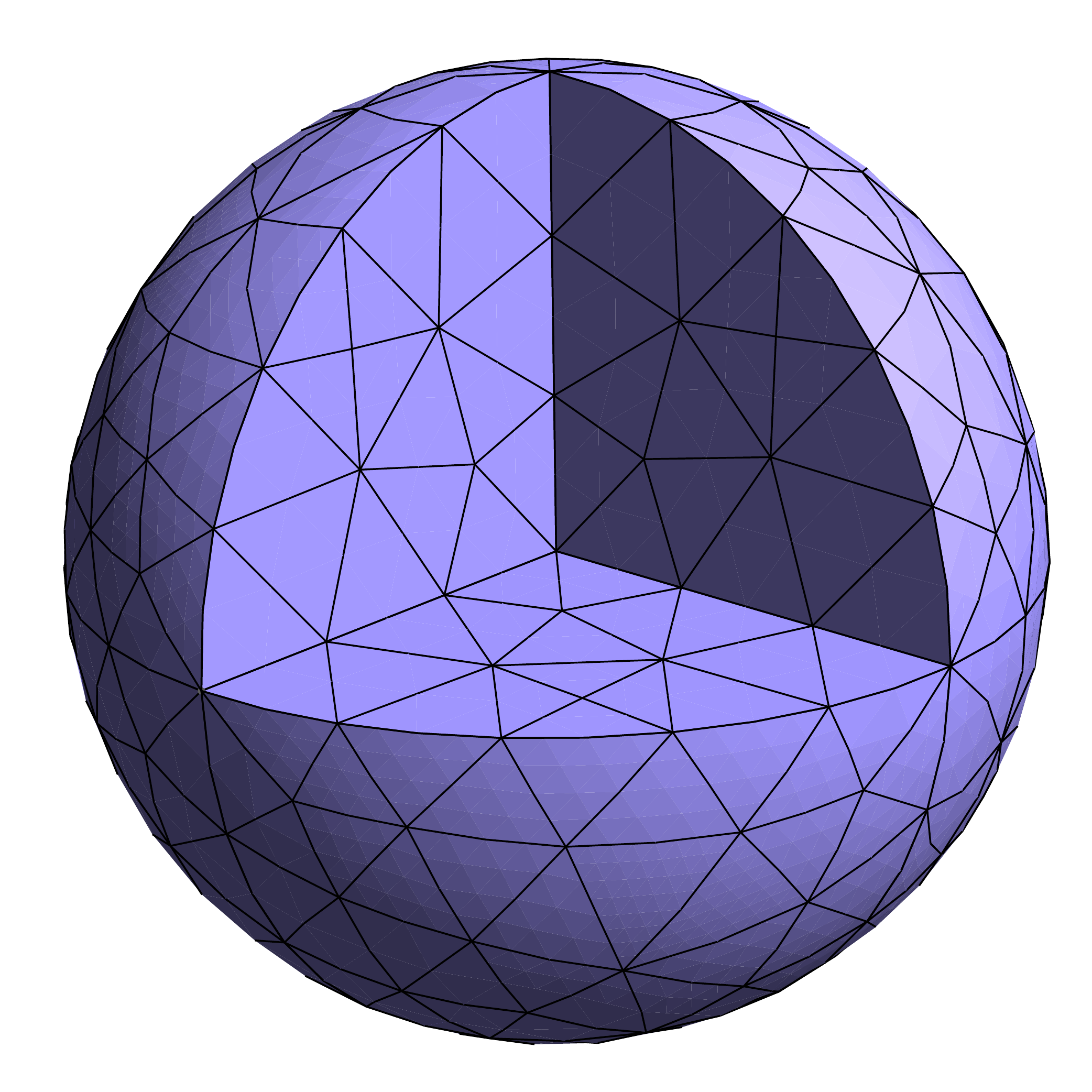}
  \caption{The cat's eye geometry}
  \label{fig:catseye}
\end{figure}

The accuracy and convergence of the integration method are tested by
implementing it in a BEM code~\cite{bem3d} and solving the Laplace
equation on a cat's eye geometry, shown in
Figure~\ref{fig:catseye}. This is a unit sphere with one octant
removed, recommended as a more stringent test of BEM codes than a
simple sphere~\cite{marburg-amini05}, since it contains
discontinuities in the geometry, as would be found, for example, in
aerodynamic calculations. The surface was meshed using
GMSH~\cite{geuzaine-remacle09}, changing the discretization length to
produce panels of varying sizes. A second mesh was produced for
comparison by splitting each second order element into six triangles,
giving a mesh of planar elements based on the same nodes. The polar
quadrature rule was selected for $\sigma<1$, with
$N_{\theta}=N_{r}=8$, and a twenty-five point symmetric
rule~\cite{wandzura-xiao03} for $\sigma\geq1$.


\begin{figure}
  \centering
  \includegraphics{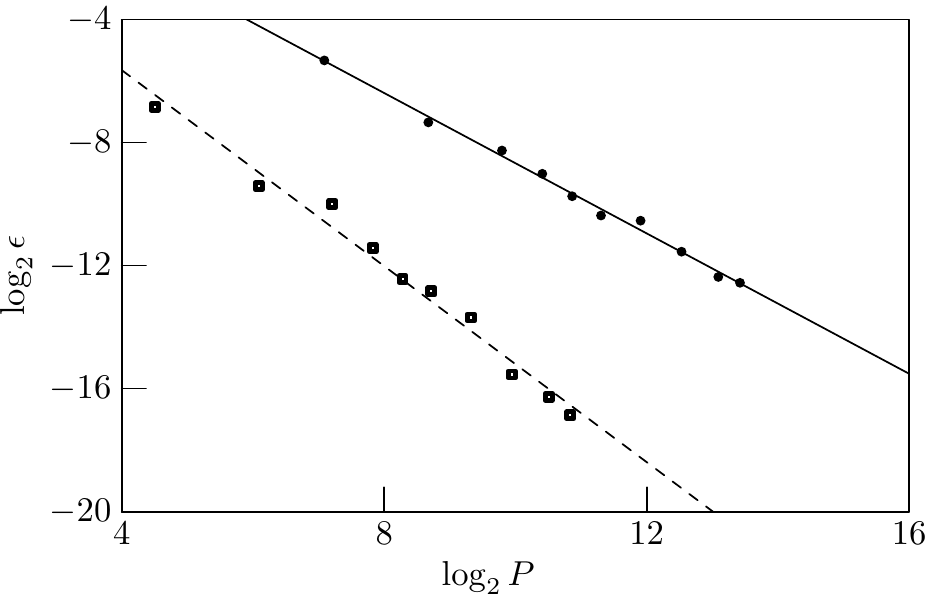} \\
  \includegraphics{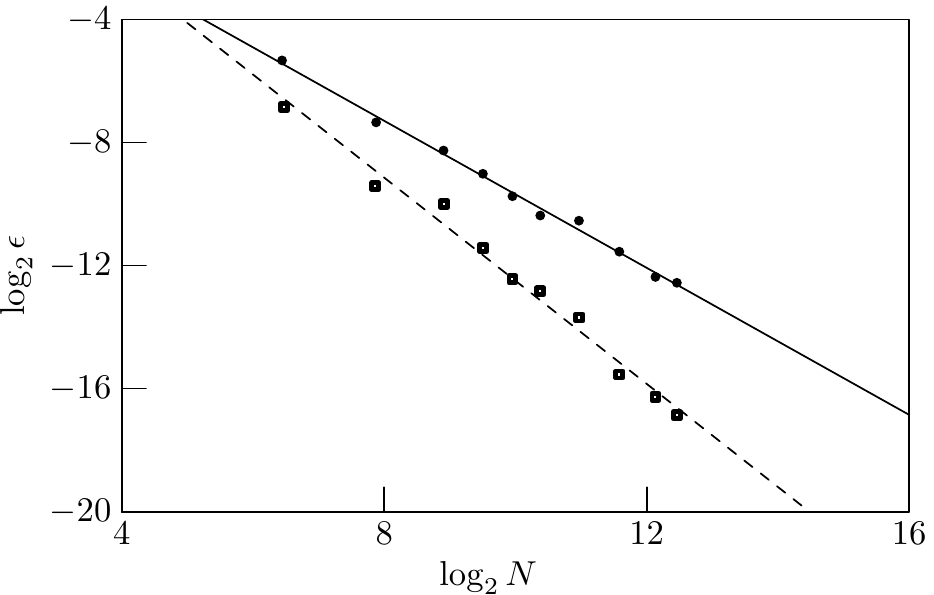}  
  \caption{Error $\epsilon$ in solution of Laplace equation on cat's
    eye against panel number $P$ (upper) and node number $N$
    (lower). Dots: linear elements; squares: quadratic elements. Upper
    plot: solid line $6.6P^{-1.1}$; dashed line: $3.9P^{-1.6}$. Lower
    plot: solid line $4.8N^{-1.2}$; dashed line: $19.3N^{-1.6}$.}
  \label{fig:laplace:error}
\end{figure}

A Neumann boundary condition was generated using a point potential
source positioned inside the surface at $(-0.2,-0.2,-0.2)$. Solving
Equation~\ref{equ:potential} for the surface potential $\phi$ gave a
result which could be compared to the result found analytically for
the point source. The error estimate is the r.m.s.\ difference between
the computed and the analytically specified data. The error is plotted
in Figure~\ref{fig:laplace:error} as a function of element number,
and, for convenience, of node number. The error is fitted using a
power law, to estimate the convergence rate, and the superior
numerical performance of the second order method is clear. Plotting
against panel number shows a similar convergence rate as in other
work~\cite[Figure~12]{willis-peraire-white06}, $P^{-1}$ for linear
panels, and $P^{-3/2}$ for quadratic. Plotting against node number,
which can be taken as a proxy for the memory requirement for the
matrix used to solve the problem, and which shows error for different
element types on the same point distribution, shows similar trends.

\section{CONCLUSIONS}
\label{sec:conclusions}
\vspace{-2pt}

A quadrature technique for second order triangular elements has been
presented and tested on a realistic geometry. It has been found that
the method is accurate and convergent, giving an error which scales as
$P^{-1.6}$ in the example tested. It is concluded that the technique
is readily implemented and can be used as a direct replacement for
existing quadratures.

\bibliography{abbrev,maths,misc,vortex,scattering,aerodynamics}

\end{document}